# Evaluation of binomial series with harmonic numbers

Khristo N. Boyadzhiev
Ohio Northern University
Department of Mathematics
Ada, Ohio 45810, USA
k-boyadzhiev@onu.edu

**Abstract**. We define a special function related to the digamma function and use it to evaluate in closed form various series involving binomial coefficients and harmonic numbers.

**MSC 2020:** 11B65; 30B10; 33B15; 40A05; 65B10

**Key words:** Digamma function, Lerch transcendent, binomial coefficient, central binomial coefficient, harmonic number, Stirling number of the first kind.

## 1. Introduction and motivation

Let $\binom{2n}{n}$ be the central binomial coefficients and let $H_n$ be the harmonic numbers, $n = 0, 1, 2, \ldots$ Studies of series with binomials and harmonic numbers have a long history (see [14, 15]). The generating function for the central binomial coefficients and related results were discussed by Lehmer in [9]. In [2] the present author computed in explicit form the generating functions for the numbers $\binom{2n}{n} H_n$ and $C_n H_n$ where $C_n = \frac{1}{n+1}\binom{2n}{n}$ are the Catalan numbers. These studies were continued further by various authors – see, for example, [5, 6, 11, 12]. Motivated by these developments, in the present paper we compute in explicit form the generating functions for the products $\binom{z}{n} H_n$ and $\binom{n+z}{n} H_n$, where $z$ can be any complex number. In particular, when $z = -1/2$ we have $\binom{-1/2}{n} = \frac{(-1)^n}{4^n}\binom{2n}{n}$, so this is a desirable generalization. For our purpose we define a special function of two variables $F(z, x)$ which can be expressed in terms of the digamma function $\psi(z) = \frac{d}{dz}\log\Gamma(z)$ and the Lerch transcendent $\Phi(x, z, 1)$. We prove that

$$\sum_{n=1}^{\infty}\binom{z}{n}H_n x^n = (x+1)^z F\left(z, \frac{x}{x+1}\right)$$

$$\sum_{n=1}^{\infty}\binom{n+z}{n}H_n x^n = (1-x)^{-z-1}\left[F(z,x) - \log(1-x)\right].$$

This makes it possible to evaluate in closed form various binomial series with harmonic numbers for different values of $z$ and $x$.

In the next section we introduce the function $F(z,x)$ and focus on its properties. We prove an integral representation for it and compute the Taylor coefficients of $F(z,x)$ with respect to the variable $z$. An important property is its functional equation given in Theorem 5.

In Section 3 we provide a number of applications given as corollaries from the main theorems 3 and 7. The new function $F(z,x)$ has interesting properties and we believe that it will be useful in other situations too.

## 2. Definition and properties of the function $F(z,x)$

We define the function of two variables

$$F(z,x) = \sum_{n=1}^{\infty}\binom{z}{n}\frac{(-1)^{n-1}x^n}{n}$$

for $z \in \mathbb{C}$ and $|x| < 1$. Clearly, $F(z,0) = 0$ and $F(0,x) = 0$. As usual $\binom{z}{n} = \frac{z(z-1)\ldots(z-n+1)}{n!}$.

When $z = p$ is a positive integer, the sum is finite and we have (see [2, entry (9.26)])

$$F(p,x) = \sum_{k=1}^{p}\binom{p}{k}\frac{(-1)^{k-1}x^k}{k} = H_p - \sum_{k=1}^{p}\frac{(1-x)^k}{k}.$$

When $z$ is not a positive integer we have the asymptotic for large $n$

$$\binom{z}{n} \approx \frac{(-1)^n}{\Gamma(-z)\, n^{z+1}}$$

and we conclude that the function $F(z,x)$ is analytic in $z$ for $|x| < 1$.

Here are some obvious properties of this function. For $z = -1$ we have $\binom{-1}{n} = (-1)^n$ and

$$F(-1, x) = -\sum_{n=1}^{\infty} \frac{x^n}{n} = \log(1-x).$$ Thus $F(-1, -1) = \log 2$. It is also clear that $F(-1, 1)$ is not defined.

However, for $\operatorname{Re} z > -1$ the above asymptotic shows that $F(z, 1)$ is well defined.

Let again $p$ be a positive integer. Then

$$\binom{-p-1}{n} = (-1)^n \binom{n+p}{n} = (-1)^n \frac{(n+p)!}{p!n!} = (-1)^n \frac{(n+1)(n+2)\dots(n+p)}{p!}$$

so clearly for $x = \pm 1$, $F(-p-1, -1)$ and $F(-p-1, 1)$ are not defined.

Now we consider further properties of $F(z, x)$.

**Lemma 1.** *For real* $-1 < x < 1$ $F(z, x)$ *has the integral representations*

$$F(z, x) = \int_0^x \frac{1-(1-t)^z}{t} dt = \int_0^1 \frac{1-(1-xt)^z}{t} dt = \int_{1-x}^1 \frac{1-u^z}{1-u} du.$$

*Proof.*

$$\sum_{n=1}^{\infty} \binom{z}{n} \frac{(-1)^{n-1} x^n}{n} = \sum_{n=1}^{\infty} \binom{z}{n} (-1)^{n-1} \int_0^x t^{n-1} dt = -\int_0^x \frac{1}{t} \left( \sum_{n=0}^{\infty} \binom{z}{n} (-1)^n t^n - 1 \right) dt$$

$$= -\int_0^x \frac{(1-t)^z - 1}{t} dt$$

For the second integral we use the substitution $t \to xt$ and for the third we use $1-t = u$.

**Note.** The digamma function has the important representation by partial fractions

$$\psi(z+1) + \gamma = \sum_{n=1}^{\infty} \left( \frac{1}{n} - \frac{1}{z+n} \right)$$

through which $\psi(z+1)$ extends by analytic continuation to all complex numbers $z \neq -1, -2, \dots$.

**Proposition 2**. *For* $\operatorname{Re} z > -1$ *and* $x = 1$ *the function* $F(z, 1)$ *is defined and we have*

$$F(z,1) = \sum_{n=1}^{\infty} \binom{z}{n} \frac{(-1)^{n-1}}{n} = \psi(z+1) + \gamma$$

(*$\gamma$ being Euler's constant*). *Therefore, $F(z,1)$ extends by analytic continuation for all $z \neq -1, -2, \ldots$ by means of the equation*

$$F(z,1) = \psi(z+1) + \gamma.$$

The proof follows immediately from Lemma 1 by setting $x \to 1$ and the well-known integral representation

$$\psi(z+1) + \gamma = \int_0^1 \frac{1-u^z}{1-u} du.$$

Note that when $z = p$ (a nonnegative integer)

$$\psi(p+1) + \gamma = \sum_{n=1}^{p} \binom{p}{n} \frac{(-1)^{n-1}}{n} = H_p.$$

In the following theorem we show that $F(z,x)$ can be expressed in terms of the digamma function and the Lerch transcendent (also called Lerch function, or Hurwitz-Lerch zeta function)

$$\Phi(t,s,a) = \sum_{k=0}^{\infty} \frac{t^k}{(k+a)^s} \quad (|t| < 1, \, a > 0)$$

(see, for example, [13]). When $s$ is an integer, $\Phi(t,s,a)$ is defined for all $a \neq 0, -1, -2, \ldots$ For our purpose we will need only the case $s = 1$, that is, $\Phi(t,1,a) = \sum_{n=0}^{\infty} \frac{t^n}{n+a}$.

**Theorem 3.** *For $z \neq -1, -2, \ldots$ and $0 < x \leq 1$ we have the representation*

$$F(z,x) = \psi(z+1) + \gamma + \log x + \sum_{n=1}^{\infty} \frac{(1-x)^{z+n}}{z+n}, \text{ or}$$

$$F(z,x) = \psi(z+1) + \gamma + \ln x + (1-x)^{z+1} \Phi(1-x, 1, z+1)$$

*For $-1 < x < 0$ and $z$ not a positive integer we have the representation*

$$F(z,x) = \psi(1-z) + \gamma + \ln(-x) + \frac{1-(1-x)^z}{z} + (1-x)^{z-1}\Phi\left(\frac{1}{1-x}, 1, 1-z\right).$$

The proof of the theorem will be given in Section 4. These representations are important for the evaluation of certain binomial series in closed form.

At this point we recall the analytic definition of the Stirling numbers of the first kind $s(n,k)$

$$\binom{z}{n} = \frac{z(z-1)...(z-n+1)}{n!} = \frac{1}{n!}\sum_{k=0}^{n} s(n,k) z^k$$

i.e. they are the coefficients of the polynomial $z(z-1)...(z-n+1)$ (see [1] and [7, p. 212]).

The next proposition gives the explicit Taylor series of $F(z,x)$ with respect to the variable $z$ centered at the origin.

**Proposition 4.** *For $|x|<1$ we have*

$$F(z,x) = \sum_{n=1}^{\infty} \frac{(-1)^{n-1} x^n}{n!n}\left\{\sum_{k=1}^{n} s(n,k) z^k\right\} = \sum_{k=1}^{\infty} (-1)^{k-1} z^k \left\{\sum_{n=1}^{\infty} \frac{(-1)^{n-k} s(n,k)}{n!n} x^n\right\}$$

$$F(z,x) = \sum_{k=1}^{\infty} (-1)^{k-1} z^k \left\{\zeta(k+1) + \sum_{j=0}^{k} \frac{(-1)^{j-1}}{j!} \operatorname{Li}_{k-j+1}(1-x) \ln^j(1-x)\right\}.$$

*Also, for $x=1$ and $|z|<1$ we have the well-known Taylor series*

$$F(z,1) = \sum_{k=1}^{\infty} (-1)^{k-1} z^k \zeta(k+1) = \psi(z+1) + \gamma.$$

*Proof.* The first equation comes directly from the definition of $s(n,k)$ after changing the order of summation. The second equation results from the first one in view of the identity

$$\sum_{n=1}^{\infty} \frac{(-1)^{n-k} s(n,k) x^n}{n!n} = \zeta(k+1) + \sum_{j=0}^{k} \frac{(-1)^{j-1}}{j!} \operatorname{Li}_{k-j+1}(1-x) \ln^j(1-x)$$

which can be found in Adamchik's paper [1, (22)] (Adamchik contributes this result to Kolbig [10]). Note that in his paper Adamchik uses the unsigned Stirling numbers of the first kind

$$\begin{bmatrix} n \\ k \end{bmatrix} = (-1)^{n-k} s(n,k).$$

With $x = 1$ in the above identity we come to the popular representation of the zeta function

$$\zeta(k+1) = \sum_{n=1}^{\infty} \frac{(-1)^{n-k} s(n,k)}{n!\, n}.$$

From the formula for the Taylor coefficients and the integral representation of $F(z,x)$ we conclude that (cf. [1, (21)])

$$\frac{(-1)^{k-1}}{k!} \left( \frac{\partial}{\partial z} \right)^k F(z,x) \bigg|_{z=0} = \frac{(-1)^k}{k!} \int_0^x \frac{\ln^k(1-t)}{t}\, dt$$

$$= \zeta(k+1) + \sum_{j=0}^{k} \frac{(-1)^{j-1}}{j!} \operatorname{Li}_{k-j+1}(1-x) \ln^j(1-x).$$

We will now prove an important property, the functional equation for $F(z,x)$. This is the second equation in the following theorem.

**Theorem 5.** *For all $z$ and all $-1 \leq x < \frac{1}{2}$ we have the representations*

$$F(z,x) = \log(1-x) - \sum_{n=1}^{\infty} \left( \frac{-x}{1-x} \right)^n \frac{1}{n} \binom{n+z}{n}$$

$$F(z,x) = \log(1-x) + F\left( -z-1, \frac{-x}{1-x} \right)$$

*or, with $-\frac{1}{2} < x \leq 1$*

$$F\left( z, \frac{x}{1+x} \right) = -\log(1+x) + F(-z-1, -x).$$

*Proof.* For the proof we use the special series identity [4, Proposition 3]

$$a_0 \log(1+x) + \sum_{n=1}^{\infty} \frac{x^n}{n} a_n = \sum_{n=1}^{\infty} \left(\frac{x}{1+x}\right)^n \frac{1}{n} \left\{\sum_{k=0}^{n} \binom{n}{k} a_k\right\}.$$

where $f(t) = a_0 + a_1 t + a_2 t^2 + \ldots$ is a function analytic in a neighborhood of the origin. We set $a_n = \binom{z}{n}$, so that $f(t) = (1+t)^z$ and $a_0 = 1$. The above identity provides the equation

$$\log(1+x) + \sum_{n=1}^{\infty} \binom{z}{n} \frac{x^n}{n} = \sum_{n=1}^{\infty} \left(\frac{x}{x+1}\right)^n \frac{1}{n} \left\{\sum_{k=0}^{n} \binom{n}{k}\binom{z}{k}\right\}$$

$$= \sum_{n=1}^{\infty} \left(\frac{x}{x+1}\right)^n \frac{1}{n} \binom{n+z}{n}$$

where the last equality uses Vandermonde's identity

$$\sum_{k=0}^{n} \binom{n}{k}\binom{z}{k} = \binom{n+z}{n}.$$

From here, replacing $x$ by $-x$ we obtain

$$\sum_{n=1}^{\infty} \binom{z}{n} \frac{(-1)^{n-1} x^n}{n} = \log(1-x) - \sum_{n=1}^{\infty} \left(\frac{-x}{1-x}\right)^n \frac{1}{n} \binom{n+z}{n}$$

which is the first equation in the theorem. The convergence of the series on the right hand side is assured by the asymptotic for large $n$

$$\binom{n+z}{n} \approx \frac{e^{z(H_n - \gamma)}}{\Gamma(z+1)}.$$

Next we use the fact that $\binom{n+z}{n} = (-1)^n \binom{-z-1}{n}$ to obtain the second equation. The third equation follows from the second. The proof is completed.

**Corollary 6.** *For* $-1 < x < \frac{1}{2}$

$$\sum_{n=1}^{\infty}\binom{n+z}{n}\frac{x^n}{n} = -\ln(1-x) - F\left(z, \frac{-x}{1-x}\right)$$

and for any positive integer $z = p$ and $|x| < 1$

$$\sum_{n=1}^{\infty}\binom{n+p}{n}\frac{x^n}{n} = -H_p - \ln(1-x) + \sum_{k=1}^{p}\frac{1}{k(1-x)^k}.$$

*Proof.* The first equation comes from the first equation of Theorem 5 by the substitution $\frac{-x}{1-x} \to x$. For the second equation we write

$$\sum_{n=1}^{\infty}\binom{n+p}{n}\frac{x^n}{n} = -\ln(1-x) - F\left(p, \frac{-x}{1-x}\right) = -\log(1-x) - H_p + \sum_{k=1}^{p}\frac{1}{k(1-x)^k}$$

by using the definition of $F(z,x)$ and the identity [3, (9.27)]

$$\sum_{k=1}^{p}\binom{p}{k}\frac{y^k}{k} = -H_p + \sum_{k=1}^{p}\frac{(1+y)^k}{k}.$$

End of proof.

Setting $x \to -1$ we come to the equation

$$\sum_{n=1}^{\infty}\binom{n+p}{n}\frac{(-1)^n}{n} = -H_p - \log 2 + \sum_{k=1}^{p}\frac{1}{k 2^k}$$

where the series on the left hand side exists in terms of Abel summability.

**Note.** In view of the expansion

$$\frac{1}{(1-x)^{p+1}} = \sum_{n=0}^{\infty}\binom{n+p}{n}x^n \quad (|x|<1)$$

we can write

$$\int_0^x \frac{1}{t}\left(\frac{1}{(1-t)^{p+1}} - 1\right)dt = \sum_{n=1}^{\infty}\binom{n+p}{n}\frac{x^n}{n}$$

and the above corollary gives the evaluation of this integral.

## 3. Evaluation of binomial series with harmonic numbers

**Theorem 7.** *For all $z$ and all real $-\dfrac{1}{2} < x < 1$ the following evaluation is true*

$$\sum_{n=1}^{\infty}\binom{z}{n}H_n x^n = (x+1)^z F\left(z, \frac{x}{x+1}\right).$$

*When $0 < x < 1$ and $z \neq -1, -2, \ldots$*

$$\sum_{n=1}^{\infty}\binom{z}{n}H_n x^n = (x+1)^z\left(\psi(z+1) + \gamma + \log\frac{x}{1+x} + \sum_{n=1}^{\infty}\frac{1}{(1+x)^{n+z}(n+z)}\right).$$

*Also, for $|x| < 1$*

$$\sum_{n=1}^{\infty}\binom{n+z}{n}H_n x^n = (1-x)^{-z-1}\left[F(z,x) - \log(1-x)\right].$$

Note that for $-\dfrac{1}{2} < x$ we have $-1 < \dfrac{x}{x+1} < 1$, so the right hand side in the first equation is well defined.

*Proof.* The proof is based on the series transformation formula (see [4, Proposition 1])

$$\sum_{n=0}^{\infty}\binom{z}{n}(-1)^n a_n x^n = (x+1)^z \sum_{n=0}^{\infty}\left(\frac{x}{x+1}\right)^n \binom{z}{n}(-1)^n \left\{\sum_{k=0}^{n}\binom{n}{k}a_k\right\}$$

where $a_n$ are the coefficients of some power series with nonzero radius of convergence. We select $a_n = (-1)^{n-1}H_n$, $a_0 = 0$. This gives the series identity

$$\sum_{n=1}^{\infty}\binom{z}{n}H_n x^n = (x+1)^z \sum_{n=1}^{\infty}\left(\frac{x}{x+1}\right)^n \binom{z}{n}\frac{(-1)^{n-1}}{n}$$

which proves the first equation in the theorem. For the second equation we use Theorem 3.

For the third equation we replace $z$ by $-z-1$ so that $\binom{n+z}{n} = (-1)^n \binom{-z-1}{n}$, apply the first equation with $-x$ in the place of $x$, and then reach for the functional equation from Theorem 5. This completes the proof.

With $z = -\frac{1}{2}$ we have $\binom{-1/2}{n} = \frac{(-1)^n}{4^n}\binom{2n}{n}$ and we obtain a new proof of Theorem 3 in [2].

**Corollary 8.** For $-1 < x \leq 1$

$$\sum_{n=0}^{\infty} \binom{2n}{n} \frac{(-1)^n H_n x^n}{4^n} = \frac{2}{\sqrt{1+x}} \log \frac{1+\sqrt{1+x}}{2\sqrt{1+x}}.$$

For the proof we take first $0 < x < 1$ and compute the right hand side $(x+1)^{-\frac{1}{2}} F\left(-\frac{1}{2}, \frac{x}{x+1}\right)$ by using the facts that $\psi\left(\frac{1}{2}\right) + \gamma = -2\log 2$, $\log \frac{x}{1+x} = \log x - 2\log\sqrt{1+x}$, and

$$2\sum_{n=1}^{\infty} \frac{1}{(1+x)^{n-\frac{1}{2}}(2n-1)} = 2\sum_{n=1}^{\infty} \left(\frac{1}{\sqrt{1+x}}\right)^{2n-1} \frac{1}{2n-1} = \log \frac{\sqrt{1+x}+1}{\sqrt{1+x}-1}$$

$$= \log \frac{(\sqrt{1+x}+1)^2}{x} = 2\log(\sqrt{1+x}+1) - \log x.$$

Putting these pieces together we come to the right hand side in the corollary. At the end we relax the restriction on $x$ to allow $-1 < x \leq 1$.

**Corollary 9.** *For every integer $z = p \geq 0$ and $|x| < 1$*

$$\sum_{n=1}^{\infty} \binom{n+p}{n} H_n x^n = \frac{1}{(1-x)^{p+1}}\left(H_p - \ln(1-x) - \sum_{k=1}^{p} \frac{(1-x)^k}{k}\right)$$

*(when $p = 0$ the sum on the right hand side is missing).*

This follows from the third equation of Theorem 7 by using the definition of $F(p, x)$ just like in the proof of Corollary 6.

**Corollary 10**. *For* $-1 \leq y \leq \frac{1}{2}$ *the function* $F(z, y)$ *has the representation*

$$F(z, y) = (1-y)^z \sum_{n=1}^{\infty} \binom{z}{n} H_n \left(\frac{y}{1-y}\right)^n = \sum_{n=1}^{\infty} \binom{z}{n} H_n y^n (1-y)^{z-n}.$$

This results from the above theorem with the substitution $y = \frac{x}{x+1}$ in the first equation.

**Corollary 11**. *For any* $\operatorname{Re} z > -1$ *we have*

$$\sum_{n=1}^{\infty} \binom{z}{n} H_n = 2^z F\left(z, \frac{1}{2}\right) = 2^z (\psi(z+1) + \gamma - \ln 2) + \frac{1}{2} \Phi\left(\frac{1}{2}, 1, z+1\right).$$

For the proof we set $x \to 1$ in Theorem 7 and also use Theorem 3. The series converges because

$$\left|\binom{z}{n}\right| \leq \frac{M}{|\Gamma(-z)| n^{\operatorname{Re} z+1}} \text{ for } z \neq 0, 1, 2, \ldots \text{ (otherwise the sum is finite).}$$

With $x = -\frac{1}{2}$ in Theorem 7 we come to the next corollary.

**Corollary 12**. *For any* $z$ *with* $\operatorname{Re} z > -1$

$$\sum_{n=0}^{\infty} \binom{z}{n} \frac{(-1)^{n-1} H_n}{2^n} = \frac{1}{2^z} \sum_{n=1}^{\infty} \binom{z}{n} \frac{1}{n} = -\frac{1}{2^z} F(z, -1).$$

*and for* $z \neq 0, 1, 2, \ldots$

$$\sum_{n=0}^{\infty} \binom{z}{n} \frac{(-1)^{n-1} H_n}{2^n} = \frac{1}{2^z} \left(\frac{2^z - 1}{z} - \psi(1-z) - \gamma - 2^{z-1} \Phi\left(\frac{1}{2}, 1, 1-z\right)\right).$$

With $x = \frac{1}{2}$ we obtain also:

**Corollary 13**. *For any* $z$

$$\sum_{n=0}^{\infty} \binom{z}{n} \frac{H_n}{2^n} = \left(\frac{3}{2}\right)^z F\left(z, \frac{1}{3}\right)$$

$$= \left(\frac{3}{2}\right)^z \left(\psi(z+1) + \gamma - \ln 3 + \left(\frac{2}{3}\right)^{z+1} \Phi\left(\frac{2}{3}, 1, z+1\right)\right).$$

By choosing $z = -\frac{1}{2}$ in Theorem 7 we obtained Corollary 8. Choosing $z = \pm\frac{1}{3}, \pm\frac{1}{4}, \pm\frac{1}{5}$ etc. we can generate other power series with harmonic numbers which could be of interest. We will demonstrate one such series in the next corollary.

**Corollary 14**. *For* $-1 < x \leq 1$

$$\sum_{n=1}^{\infty} (-1)^{n-1} \Gamma\left(n - \frac{1}{4}\right) H_n \frac{x^n}{n!}$$

$$= 4\Gamma\left(\frac{3}{4}\right) \sqrt[4]{x+1} \left(4 - \frac{\pi}{2} - 3\log 2 - \frac{4}{\sqrt[4]{x+1}} + \log\frac{x(\sqrt[4]{x+1}+1)}{(x+1)(\sqrt[4]{x+1}-1)} + 2\arctan\frac{1}{\sqrt[4]{x+1}}\right).$$

*Proof.* We start the proof with $0 < x < 1$ and at the end we drop this restriction. We use Theorem 7 with $z = \frac{1}{4}$ in the second equation of Theorem 7. We also use the facts

$$\binom{1/4}{n} = \frac{(-1)^{n-1}}{4n!\Gamma(3/4)} \Gamma\left(n - \frac{1}{4}\right), \quad \psi\left(1 + \frac{1}{4}\right) + \gamma = 4 - \frac{\pi}{2} - 3\log 2$$

$$\sum_{n=1}^{\infty} \frac{t^{4n+1}}{4n+1} = -t + \frac{1}{4}\log\frac{1+t}{1-t} + \frac{1}{2}\arctan t \quad (|t| < 1).$$

A reference for the above series is Hansen's table [8, entry 5.10,20]. In our case $t = (x+1)^{-1/4}$. After computing the expression on the right hand side, we can easily see that it is defined for $-1 < x \leq 1$ as

$$\lim_{x \to 0} \log \frac{x(\sqrt[4]{x+1}+1)}{(x+1)(\sqrt[4]{x+1}-1)} = 3\log 2.$$

With this the proof is completed.

## 4. Proof of Theorem 3

For the proof we use the integral representation from Lemma 1.

$$F(z,x) = \int_{1-x}^{1} \frac{1-u^z}{1-u} du$$

When $0 < x \leq 1$ we expand $(1-u)^{-1}$ in geometric series and then integrate between $1-x$ and $1$

$$F(z,x) = \sum_{n=1}^{\infty} \int_{1-x}^{1} (1-u^z) u^{n-1} du = \sum_{n=1}^{\infty} \left\{ \frac{1-(1-x)^n}{n} - \frac{1-(1-x)^{z+n}}{z+n} \right\}$$

$$= \sum_{n=1}^{\infty} \left( \frac{1}{n} - \frac{1}{z+n} \right) - \sum_{n=1}^{\infty} \frac{(1-x)^n}{n} + \sum_{n=1}^{\infty} \frac{(1-x)^{z+n}}{z+n}.$$

The first sum is $\psi(z+1) + \gamma$ and we obtain the first equation in the proposition. The second equation follows from the definition of the Lerch transcendent since replacing $n$ by $n+1$ we can write

$$\sum_{n=1}^{\infty} \frac{(1-x)^{z+n}}{z+n} = (1-x)^z \sum_{n=1}^{\infty} \frac{(1-x)^n}{z+n} = (1-x)^{z+1} \sum_{n=0}^{\infty} \frac{(1-x)^n}{n+1+z}.$$

The proof of the first part is completed.

For the second part (when $-1 \leq x < 0$) we reason differently, expanding $(1-1/u)^{-1}$ in geometric series, as this time $u \geq 1$

$$F(z,x) = \int_{1-x}^{1} \frac{1-u^z}{1-u} du = \int_{1-x}^{1} \frac{u^{z-1}-1/u}{1-1/u} du = \sum_{n=0}^{\infty} \int_{1-x}^{1} (u^{z-1} - u^{-1}) u^{-n} du$$

$$= \sum_{n=0}^{\infty} \int_{1-x}^{1} (u^{z-n-1} - u^{-n-1}) du = \int_{1-x}^{1} (u^{z-1} - u^{-1}) du + \sum_{n=1}^{\infty} \int_{1-x}^{1} (u^{z-n-1} - u^{-n-1}) du$$

$$= \frac{1-(1-x)^z}{z} + \log(1-x) + \sum_{n=1}^{\infty} \left\{ \frac{1-(1-x)^{z-n}}{z-n} + \frac{1-(1-x)^{-n}}{n} \right\}$$

$$= \frac{1-(1-x)^z}{z} + \log(1-x) + \sum_{n=1}^{\infty} \left\{ \frac{1}{n} - \frac{1}{n-z} \right\} + \sum_{n=1}^{\infty} \frac{(1-x)^z}{(1-x)^n(n-z)} - \sum_{n=1}^{\infty} \frac{1}{(1-x)^n n}.$$

Here

$$\sum_{n=1}^{\infty}\frac{(1-x)^z}{(1-x)^n(n-z)} - \sum_{n=1}^{\infty}\frac{1}{(1-x)^n n} = (1-x)^{z-1}\sum_{n=0}^{\infty}\left(\frac{1}{1-x}\right)^n \frac{1}{n+1-z} + \log\left(1-\frac{1}{1-x}\right)$$

$$= (1-x)^{z-1}\Phi\left(\frac{1}{1-x}, 1, 1-z\right) + \log\left(\frac{-x}{1-x}\right).$$

The logarithms come together and we obtain the last equation of the theorem. The proof is completed.